\magnification=1200
\baselineskip=16pt
\voffset=-1.7cm
\vsize=19.6cm
\overfullrule=0pt
\def\I{I\!\!P}
\font\c=cmr8
\font\goth=eusm10
\def\fake{\hbox{\goth I}}
\def\medtype{\let\rm=\tenrm 
\baselineskip=12pt minus 1pt
\rm}
\nopagenumbers
\headline={\ifnum\pageno=1 \else{\ifodd\pageno\rightheadline 
\else\leftheadline\folio\fi}\fi}
\def\rightheadline{\c\hfil FOCAL LOCI OF FAMILIES AND THE GENUS 
OF CURVES ON SURFACES\hfil\folio}  
\def\leftheadline{\c\hfil CHIANTINI-LOPEZ\hfil}
\voffset=1\baselineskip
\font\ms=msbm10

\font\msq=msbm7

\def\qed{\vrule height 1.1ex width 1.0ex depth -.1ex }
\def\mapright#1{\smash{\mathop{\longrightarrow}\limits^{#1}}}

\centerline{\bf FOCAL LOCI OF FAMILIES}
\centerline{\bf AND THE GENUS OF CURVES ON SURFACES} 
\vskip .6cm

\centerline{LUCA CHIANTINI\footnote*{Research  
partially supported by the  MURST national project ``Geometria Algebrica"; 
the authors are members of GNSAGA of CNR.} AND 
ANGELO FELICE LOPEZ*\footnote\null{{\it
\noindent 1991  Mathematics Subject Classification:} Primary 14J29. 
Secondary 32H20, 14C20.}} 

\vskip .4cm \midinsert \narrower
\hsize=13.3cm
\medtype{

\noindent ABSTRACT: In this article we apply the classical method 
of focal loci of families to give a lower bound for the genus of 
curves lying on general surfaces. First we translate and reprove
Xu's result that any curve $C$ on a general surface in $\I^3$ of 
degree $d \geq 5$ has geometric genus $g > 1 + 
\hbox{deg} C (d - 5) / 2$. Then we prove a similar lower bound for 
the curves lying on a general surface in a given component of the 
Noether-Lefschetz locus in $\I^3$ and on a general projectively 
Cohen-Macaulay surface in $\I^4$.}
\endinsert

\vskip .9cm
\magnification=1200
\baselineskip=16pt
\vsize=19.6cm

\noindent {\bf 1. INTRODUCTION}

\vskip .5cm

The theory of singular curves lying in a projective variety $X$ 
has been extensively studied from the beginnings of Algebraic
Geometry;  however, even when $X$ is a smooth surface, many basic
questions still remain  open.
Recently, the interest about these arguments grew, essentially for 
two reasons: on the one hand, the theory of strings of nuclear 
physicists deals with the enumerative geometry of rational curves
contained in some projective threefolds; on the other, the study of
singular curves is  naturally related with the hyperbolic geometry
of complex projective  varieties. Let us recall briefly this last
setting. A compact complex manifold $M$ is said to be hyperbolic 
(in the sense of Kobayashi [K] or Brody [B]) if there are no
nonconstant entire holomorphic maps $f : \hbox{\ms C} \to M$. An
intriguing question that lies in the intersection between
differential and algebraic geometry is to characterize which
projective algebraic varieties $X$ over the complex field are
hyperbolic. One necessary condition, in the above case, is that
there are  no nonconstant holomorphic maps $f : A \to X$ from an
abelian variety $A$ to
$X$ and it has been conjectured by Kobayashi [K] and Lang [La] that 
this condition is in fact sufficient. 

An approach to the themes around this 
conjecture was given by Demailly in his Santa Cruz notes [D], in 
which he proposed an intermediate step:

{\bf Definition (1.1).} A smooth projective variety $X$ is 
said to be {\it algebraically hyperbolic} if there exists a real
number $\epsilon > 0$ such that every algebraic curve $C \subset X$
of geometric genus $g$ and degree $d$ satisfies $2g - 2 \geq
\epsilon d$.

Demailly proved that hyperbolic implies algebraically hyperbolic and
that the latter does not allow the existence of nonconstant
holomorphic maps from an abelian variety to $X$. Then, in view of
Kobayashi-Lang's conjecture, it  becomes relevant to check which
projective varieties are algebraically hyperbolic. Already in the
case of surfaces a complete answer appears far from reach. In
$\I^3$, the works of Brody, Green [BG], Nadel [N] and El Goul [EG]
showed that for all integers $d \geq 14$ there exist hyperbolic
surfaces of degree $d$. Very recently Demailly and El Goul [DEG]
proved that a {\it general} surface (in the countable Zariski topology)
in $\I^3$ of degree at least 42 is
hyperbolic. On the other hand, Clemens [Cl] proved that a {\it
general} surface of degree at least
$6$ in $\I^3$ is algebraically hyperbolic. Clemens' argument was 
extended by Ein ([E1], [E2]) to the case of complete intersections 
in higher dimensional varieties and recently improved and simplified
by Voisin in [V] (see also [CLR] for a very recent improvement). 
In 1994, Xu [X1] provided a more precise lower
bound for the geometric genus  of curves in any linear system over a
general surface in $\I^3$, which is in fact sharp in some cases. His
proof is based on a clever investigation of the equations defining
singular curves on surfaces moving in $\I^3$; however it is involved
in explicit hard computations with local coordinates. It was when trying 
to understand Xu's method from a global point of view, that we
realized its connection with the focal locus of a family of curves.
In general, the theory of ``focal loci" for a family of varieties
was classically developed by  C.Segre [S] and recently rephrased in
a modern language by Ciliberto and  Sernesi [CS] (see also [CC]).
These loci play an important role in  differential projective
geometry, and hence they have a quite natural  involvement in the
study of algebraic hyperbolicity of projective varieties.

Our first task, in the present article, has been to translate Xu's 
local analysis with a global property of the focal locus of a family
of curves (Proposition (2.4)). This property turns out to be simple
and powerful  enough to get interesting applications.

In the case of general surfaces in $\I^3$, we give a short proof of 
one of the main theorems of Xu [X1, Theorem 2.1], only by means of
focal loci.

\vskip .3cm

{\bf Theorem (1.2).}  {\sl On a general surface $S$ of degree $d 
\geq 5$ in $\I^3$, there are no reduced irreducible curves $C$ of
geometric genus 
$g \leq 1 + \hbox{deg} C (d - 5) / 2$. In particular, for $d \geq 
6$, $S$ is algebraically hyperbolic.}

\vskip .3cm

\noindent An interesting consequence is that we reobtain a proof of 
Harris' conjecture: a general quintic surface in $\I^3$ does not
contain rational or elliptic curves (see  also Remark (3.2)). It
should be noted that from  the articles of Ein [E1],  [E2] and
Voisin [V], it follows that on a general surface of degree $d \geq
5$ there are no reduced irreducible curves with 
$g \leq \hbox{deg} C (d - 5) / 2$.

The next possibility for surfaces in $\I^3$, which was also a 
starting point of our work, is to analyze the question of algebraic
hyperbolicity for  surfaces that are general in a given component of
the Noether-Lefschetz locus, that is the locus of smooth surfaces of
degree $d \geq 4$ in $\I^3$ whose  Picard group is {\it not}
generated by the hyperplane bundle. The problem was partly motivated
by the simple observation that a surface of general type which is
not algebraically hyperbolic and does not contain rational or
elliptic curves,  would be a counterexample to Kobayashi-Lang's
conjecture. The simplest examples of surfaces in a component of the
Noether-Lefschetz locus are those containing a fixed curve $D
\subset \I^3$. The method of focal loci proved useful also in this 
case:

\vskip .3cm

{\bf Theorem (1.3).} {\sl Let $D$ be an integral curve in $\I^3$ and
let $s, d$ be two integers such that $d \geq s + 4$ and\par
(i) there exists a surface $Y \subset 
\I^3$ of degree $s$ containing $D$,

(ii) the general element of the linear system 
$|{\cal O}_{Y}(dH - D)|$ is irreducible.  

\noindent Let $S$ be a general surface of degree $d$ in $\I^3$, 
containing $D$ and set $Y \cap S = D \cup D'$. Then $S$ contains no reduced
irreducible curves $C \neq D, D'$ of geometric genus 
$g < 1 + \hbox{deg} C (d - s - 5) / 2$. In particular, for $d \geq s + 6$ 
and $g(D), g(D') \geq 2$, $S$ is algebraically hyperbolic.}

\vskip .3cm

\noindent This result should be compared with what can be obtained 
adapting the methods of [V], namely that on $S$ there are no reduced
irreducible curves with
$g < 1+ \hbox{deg} C (d - \alpha - 4) / 2$, where $\alpha$ is a 
degree in which the homogeneous ideal of $D$ is generated [Voisin,
{\it priv. comm.}].

Finally we give another application to the case of projectively 
Cohen-Macaulay surfaces in $\I^4$, where the methods of Ein and
Voisin do not seem to apply easily, because they are not complete intersections. 
To establish the notation, let
$S \subset \I^4$ be a projectively Cohen-Macaulay surface and
consider the minimal free resolution of its ideal sheaf
$\fake_{S}$   
$$0 \to \bigoplus\limits_{i=1}^m {\cal
O}_{\I^4}(-d_{2i}) \mapright{\phi}  \bigoplus\limits_{j=1}^{m+1}
{\cal O}_{\I^4}(-d_{1j}) \mapright{\psi}  \fake_{S} \to 0$$ 
\noindent where we assume $d_{2i} \geq d_{2,i+1}, d_{1j} \geq 
d_{1,j+1}$. Now set $u_{ij} = d_{2i} - d_{1j}$;  note that the order
chosen implies $u_{i+1,j} \leq u_{ij} \leq u_{i,j+1}$. We have

\vskip .3cm

{\bf Theorem (1.4).}  {\sl On a general projectively Cohen-Macaulay 
surface $S$ in $\I^4$ such that $u_{m,m+1} \geq 6$ there are no
reduced irreducible  curves $C$ of geometric genus $g < 1 +
\hbox{deg} C (u_{m,m+1} - 7) / 2$. In particular, for $u_{m,m+1}
\geq 8$, $S$ is algebraically hyperbolic.}

\eject

\noindent {\bf 2. SOME BASIC FACTS ABOUT FOCAL SETS} 

\vskip .5cm

In this section we recall the construction of the focal set of a
projective family of curves. We refer to [CC] for more details.

Let ${\cal X} \to B$ be a family of hypersurfaces $X_{t} 
\subset \I^n, t \in B$ and let ${\cal C} \to B$ be a family of 
curves of geometric genus $g$ such that $C_{t} \subset X_{t}$, for
every $t \in B$. Let $\pi : {\cal C} \to \I^n$ be the corresponding
map and denote by $z({\cal C})$ the dimension of the image of $\pi$.
Without loss of generality we can assume,  by shrinking $B$, the
existence of a global desingularization 
$$\matrix{\sigma : \widetilde{\cal C} \to {\cal C} \cr  \ \ \ \ 
\searrow  \ \swarrow \cr  \ \ \ \ B \cr}$$\par 
\noindent of all the fibers of ${\cal C}$. Let $0 \in B$ be a 
general point and denote by $X_{0}, C_{0}, \widetilde C_{0}$ and
$\sigma_{0} : \widetilde C_{0} \to C_{0}$ the corresponding fibers.
We have the basic

{\bf Proposition (2.1).} {\sl Let $s :\widetilde{\cal C}\to {\cal C}
\hookrightarrow B \times \I^n$ be the composition, $\cal N$  the 
cokernel of the induced map $T_{\widetilde{\cal C}}\to s^*(T_{B
\times \I^n})$ (the  {\it normal sheaf} to $s$) and let $N_{0}$ be
the restriction of $\cal N$ to  the fiber $\widetilde{\cal C}_{0}$.
Then

\noindent (a) $N_{0}$ is the normal bundle to the composition 
$s_{0} = \pi \circ \sigma_{0}$, i.e. the cokernel of the map $\null
\ \ \ \ \ T_{\widetilde{{C_0}}} \to s_{0}^*(T_{\I^n})$;

\noindent (b) the family ${\cal C} \to B$ induces a
characteristic map 
$$\chi_{0} : T_{B} \otimes {\cal O}_{\widetilde C_{0}} \to N_{0}$$
\noindent where $T_{B}$ is the tangent space to $B$ at $0$; its 
rank at a general point of $\widetilde C_{0}$ is $z({\cal C})-1$.} 

\noindent {\it Proof:}  The first fact is Proposition 1.4 of [CC]; 
the second fact is shown in [CC, p.98]. \ \qed

{\bf  Definition (2.2).} The {\it global focal set} $F_{0}$ of the 
family ${\cal C} \to B$ is the locus defined on $\widetilde C_{0}$ 
by $\bigwedge^{n-1}\chi_{0} = 0$.

\noindent The global focal set $F_{0}$ has the following simple but 
useful properties.

{\bf Proposition (2.3).} {\sl Let $P \in \widetilde C_{0}$ be a 
point such that $Q = \sigma_{0}(P)$ is a smooth point of $C_{0}$ and
$X_{0}$. We have

\noindent (i) $F_{0} = \widetilde C_{0}$  if and only if 
$z({\cal C}) < n$;

\noindent (ii) if $z({\cal C}) = n$ and $Q$ is a fixed point of 
${\cal X}$, then $P \in F_{0}$.} 

\noindent {\it Proof:} (i) is a consequence of Proposition (2.1)(b).
 To see (ii) observe that $s$ factorizes via the inclusion ${\cal X}
\hookrightarrow  B \times \I^n$; hence we have a map 
$$ N_{0} \to \sigma_{0}^*(N_{X_{0} / \I^n} \otimes 
{\cal O}_{C_{0}})$$
\noindent whose restriction to $Q$ is surjective, since $Q$ is a 
smooth point of $X_{0}$. Now, if $Q$ is a fixed point of ${\cal X}$,
then $\chi_{0}(T_{B} \otimes {\cal O}_{\widetilde C_{0},P})$ is in
the rank $n-2$ kernel of the previous map. \ \qed

The following result is the interpretation of Xu's method ([X1], 
[X2], [X3]) in terms of focal sets.

{\bf Proposition (2.4).} {\sl Let $F_{s}$ be the subset of the 
global focal set $F_{0}$, supported at the points $P \in \widetilde
C_{0}$ which map to  the regular locus of $C_{0}$ and assume
$z({\cal C}) = n$. Then 
$$\hbox{deg} F_{s} \le 2g-2 + (n+1) \hbox{deg} C_{0}.$$} 
\noindent {\it Proof:} The first Chern class of $N_{0}$ is 
$2g - 2 + (n+1) \hbox{deg} C_{0}$, by construction. In general,
$N_{0}$ may  have some torsion subsheaf ${\cal T}$, but in any case
${\cal T}$ is supported at the {\it cuspidal points} of $\widetilde
C_{0}$, which map to the singular locus of
$C_{0}$ (see [CC, 1.6]). Let $N_{0}'= N_{0} / {\cal T}$ and consider
the composition  $\chi_{0}':T_{B} \otimes 
{\cal O}_{\widetilde C_{0}} \to N_{0} \to N_{0}'$. By our
assumptions on $z(\cal C)$ this map is generically surjective and
$F_{s}$ is contained in the locus where $\chi_{0}'$ drops rank. It
follows that $\hbox{deg} F_{s} \leq c_{1}(N_{0}') \leq
c_{1}(N_{0})$. \ \qed 

\vskip .5cm

\noindent {\bf 3. THE GENUS OF EFFECTIVE DIVISORS VIA THE FOCAL SET} 
 
\vskip .5cm

The strategy that we will employ to show that on a general member
of a given family of surfaces in $\I^n$, ${\cal S} \to B$, there are
no curves of ``low" geometric genus is the following. Supposing
there is a family ${\cal C} \to B$ of curves on each surface, first
find a collection of  subfamilies $B(U) \subset B$ such that their
tangent spaces span $T_{B}$ and  that there is a family of
hypersurfaces ${\cal X}(U)$ containing the curves and with many
fixed points. Then find one such $U$ so that the corresponding
family of curves ${\cal C}(U)$ satisfies $z({\cal C}(U)) = n$ and
apply Propositions (2.3) and (2.4) to get a lower bound for the
geometric genus. 

We first record an elementary lemma on vector spaces that we will 
use in all the proofs.

{\bf Lemma (3.1).} {\sl Let $f: V \to W$ be a linear map of finite 
dimensional vector spaces and suppose $\hbox{dim} f(V) > z$. Let
$\{V_{i}\}$ be a family  of subspaces of $V$, such that $V$ is
generated by $\bigcup V_{i}$ and assume that for any pair of
subspaces $V_{i}, V_{j}$ in the family, we can find a chain of
subspaces of the family $V_{i} = V_{i_{1}}, V_{i_{2}}, \dots ,
V_{i_{k}} = V_{j}$ such that for all $h$, $\hbox{dim} f(V_{i_{h}}
\cap V_{i_{h+1}}) \geq z$. Then there is a $V_{i}$ in the family
with $\hbox{dim} f(V_{i}) > z$.}

\noindent {\it Proof:} Assuming that $\hbox{dim} f(V_{i}) \leq z$ 
for all $i$ we will prove that, for any  $V_{i}, f(V_{i})$ generates
$f(V)$, a contradiction. Indeed, for any $j$, take a chain $V_{i} =
V_{i_{1}},  V_{i_{2}}, \dots , V_{i_{k}} = V_{j}$ which links
$V_{i}$ and $V_{j}$. We  have $\hbox{dim} f(V_{i}) \leq z \leq
\hbox{dim} f(V_{i} \cap V_{i_{2}}) \leq
\hbox{dim} f(V_{i_{2}}) \leq z$ and also $z \leq \hbox{dim}  
f(V_{i} \cap V_{i_{2}}) \leq \hbox{dim} f(V_{i}) \leq z$; hence
$f(V_{i}) = f(V_{i} \cap V_{i_{2}}) = f(V_{i_{2}})$. Repeating the
argument, we finally get $f(V_{i}) = f(V_{j})$. Since $\bigcup
f(V_{i})$ generate $f(V)$, we are done. \ \qed

Our first application will be to give a new proof of Xu's theorem 
([X1, Theorem 2.1]).

\noindent {\it Proof of Theorem (1.2):} Take a family ${\cal S}\to 
B$ of surfaces of degree $d$ in $\I^3$, with $B$ dense in 
$\I H^0({\cal O}_{\I^3}(d))$. Let ${\cal C}\to B$ be a family of 
reduced irreducible curves, such that for all $t \in B$, the fiber
$C_{t} \subset S_{t}$ has geometric genus $g$ and let $\widetilde
{\cal C} \to B$ be a global desingularization of the fibers. Fix
$0\in B$ general and suppose, without loss of generality, that the
authomorphisms of $\I^3$ act on $B$ and moreover that $S_{0}$ does
not contain any line. The action of PGL(3) on $B$ shows  that the
image of ${\cal C}$ in $\I^3$ cannot be contained in any fixed 
surface; hence $z({\cal C}) = 3$ and therefore the characteristic
map of
$\widetilde{\cal C}$ has rank 2 at a general point of 
$\widetilde{\cal C}$ by Proposition (2.1). For all  surfaces $U$ of
degree $d-1$, transversal to $C_{0}$, and for $P \in C_{0} - U$
general, let $B(U)$ be the subvariety of $B$ parametrizing surfaces
which contain $U \cap S_{0}$ and let $B(U,P)$ be  the subvariety of
$B(U)$ parametrizing surfaces passing through $P$. Denote by
$T(U)$ and $T(U,P)$ their tangent spaces at $0$ and by $\widetilde 
{\cal C}(U), \widetilde {\cal C}(U,P)$ the corresponding families 
of curves. Note that $\hbox{dim} B(U) \geq 3, \ \hbox{dim} B(U,P)
\geq 2$. We claim that the characteristic map of $\widetilde{\cal
C}(U,P)$ has rank 2 for $U,P$ general. To prove this, let $U,U'$ be
two monomials of degree $d-1$ which  differ only in degree one, that
is, $M = $ l.c.m.$(U,U')$ has degree $d$; then $B(U) \cap B(U')$
contains the pencil $S_{0} + \lambda M$, which defines a non-trivial
deformation of $C_{0}$, since $C_{0}$ is not in the base locus.  It
follows that any pair $T(U), T(U')$, with $U,U'$ monomials, can be
linked  by a chain of subspaces of this type, in such a way that 
the intersection of  two consecutive elements of  the chain has
non-trivial image under the characteristic map. Since $T_{B}$ is
generated by the tangent vectors to the varieties $B(U)$, with $U$
monomial, and the characteristic map of
$\widetilde{\cal C}$ has rank 2 at a general point, we get by Lemma 
(3.1) that the characteristic map of $\widetilde{\cal C}(U)$ has 
also rank 2 for some $U$. Similarly, fixing a general $U$, for any
smooth points $P,P' \in C_{0} - U$, the intersection $B(U,P) \cap
B(U,P')$ contains a pencil $S_{0} + \lambda M$ with 
$M = U \cdot$(some plane through $PP'$), which induces a non 
trivial deformation of $C_{0}$ (as $C_{0}$ is not a line). Hence 
applying Lemma (3.1) again one gets the claim. 

\noindent Now look at the focal locus of
$\widetilde{\cal C}(U,P)$: it contains the inverse image of 
$C_{0} \cap U$ and $P$, and therefore, by Propositions (2.3) and
(2.4), we have  
$$2g - 2 \geq (d - 1) \hbox{deg} C_{0} + 1 - 4 \hbox{deg} C_{0}$$  
\noindent that is the required inequality on $g$, as $\hbox{deg} 
C_{0}$ is a multiple of $d$ by the Noether-Lefschetz theorem. \ \qed

\noindent {\bf (3.2)}{\it Remark.} Note that the above proof makes 
use of the Noether-Lefschetz theorem only at the very last line, and
in fact only for $d$ even, deg$C$ odd. In particular our proof of
Harris' conjecture is independent of the Noether-Lefschetz theorem,
while Xu's proof makes essential  use of the fact that on a general
surface every curve is a complete intersection. However Xu gets $g
\geq 3$ for a curve on a general quintic.

We now consider a fixed integral curve $D \subset \I^3$ and the 
family of surfaces of degree $d \geq s + 4$ containing it, where $s$
and $d$ are as in  the hypotheses (i) and (ii) of Theorem (1.3).

\noindent {\it Proof of Theorem (1.3):} Act with PGL(3) on $D$ and 
let $\cal D$ be the corresponding family of curves. Let ${\cal S}
\to B$ be the family of surfaces of degree $d$ in $\I^3$, containing
some curve in ${\cal D}$ and let, as above, ${\cal C} \to B$ be a
family of reduced irreducible curves, such that for all $t \in B$,
the fiber $C_{t} \subset S_{t}$ has geometric genus $g$; let
$\widetilde {\cal C}\to B$ be a global desingularization of the
fibers. Fix $0\in B$ general, call $S_{0}$ and
$C_{0}$ the corresponding fibers and suppose again that PGL(3)  
acts on $B$. Let $Y_{0}$ be the surface as in (i) containing 
$D_{0}$. Suppose first that $C_{0} \not\subset Y_{0}$. The action
of  PGL(3) gives, as above, $z({\cal C}) = 3$;  hence the
characteristic map of $\widetilde{\cal C}$ has rank 2 by Proposition
(2.1). For all surfaces $U$ of degree $d - s - 1$ transversal to
$D_{0}$ and $C_{0}$, let $B(U)$ be the subvariety of
$B$ parametrizing surfaces which contain $U \cap S_{0}$; let $T(U)$
be its tangent space at $0$ and $\widetilde {\cal C}(U)$ the
corresponding family of curves. Note that $\hbox{dim} B(U) \geq 2$,
for $D_{0}$ is contained in many surfaces of degree $s + 1$. As in
the proof of Theorem (1.2) let $U,U'$ be two monomials of degree
$d-s-1$ which differ only in degree one, that is, $M = $
l.c.m.$(U,U')$ has degree $d-s$. Then $B(U) \cap B(U')$ contains the
pencil $S_{0} + \lambda Y_{0} M$, which defines a non-trivial
deformation of $C_{0}$, since $C_{0}$ is not in the base locus.
Therefore, by Lemma (3.1), we get that the characteristic map of
$\widetilde{\cal C}(U)$ has rank 2 for $U$ general. Now look at 
the focal locus of $\widetilde{\cal C}(U)$: it contains the inverse
image of $C_{0}\cap U$; hence  by Propositions (2.3) and (2.4), we
get  
$$2g - 2 \geq (d - s -1) \hbox{deg} C_{0} - 4 \hbox{deg} C_{0}.$$   
\noindent On the other hand if $C_{0} \subset Y_{0}$, then, by (ii), 
$C_{0}$ is the residue of $D_0$ in the complete intersection of $S_{0}$ and 
$Y_{0}$. We have then proved that for a  general curve $D_{0} \in {\cal D}$, 
a general surface of degree $d$  containing $D_{0}$ satisfies the assertion
of the theorem. Since all the  curves in ${\cal D}$ are projectively
isomorphic to $D$, the theorem follows. \ \qed

\noindent {\bf (3.3)}{\it Remark.} The bound of Theorem (1.3) can 
be improved in some cases by fixing some base points in the family
$B(U)$, as we did in the proof of Theorem (1.2). 
 
\noindent {\bf (3.4)}{\it Remark.} When $D$ and/or $D'$ are rational or
elliptic, we still have the algebraic hyperbolicity of the open surface $S -
D \cup D'$, for $d \geq s + 6$. More than that, we know that every
nonconstant map $\hbox{\ms C} \to S$ has image contained in $D \cup D'$.

\noindent {\bf (3.5)}{\it Remark.} The proofs of Theorems (1.2) and
 (1.3) can be easily extended both to the case of general complete
intersection surfaces $S \subset \I^r$ and to the case of general
complete intersection surfaces $S \subset \I^r$ containing a fixed
curve $D$.

Finally we deal with the case of projectively Cohen-Macaulay 
surfaces $S \subset \I^4$. Let us recall some notation from the
introduction. We denote by $M_{S} =  [A_{ij}]$ and $[F_{1}, \ldots ,
F_{m+1}]$ the matrices of the  maps $\phi$ and $\psi$, respectively,
appearing in the resolution of the ideal sheaf of $S$. As is well
known, either $A_{ij}$ is a polynomial of degree
$u_{ij}$ if $u_{ij} \geq 0$ or $A_{ij} = 0$  if $u_{ij} < 0$ and, 
by the Hilbert-Burch theorem, we can assume that $F_{j}$ is the
determinant of the minor obtained by removing the j-th column from
$M_{S}$.

\noindent {\it Proof of Theorem (1.4):} Take a family ${\cal S} \to 
B$ of projectively Cohen-Macaulay surfaces $S_{t} \subset \I^4, t
\in B$ with $B$  dense and let ${\cal X} \to B$ be the corresponding
family of hypersurfaces $X_{t} \subset \I^4$ defined by the minors
$F_{1,t}$ of the matrix $M_{S_{t}}$. Note that by [Ch] we can assume
that $\hbox{dim} Sing (X_{t}) =  0$, as $u_{ii} > 0$ by minimality.
Let ${\cal C}\to B$ a family of reduced irreducible curves, such
that for all $t \in B$, the fiber $C_{t}$ has geometric genus $g$
and is contained in $S_{t}$ and let $\widetilde {\cal C}\to B$  be a
global desingularization of the fibers. Fix $0 \in B$ general and
suppose that PGL(4) acts on B. The action of PGL(4) on $B$ shows
that the image of ${\cal C}$ in $\I^4$ cannot  be contained in any
fixed proper subvariety; hence $z({\cal C}) = 4$  and the
characteristic map of $\widetilde{\cal C}$ has rank 3 by Proposition
(2.1).  For all hypersurfaces $U$ of degree $u = u_{m,m+1} - 2$,
transversal to $C_{0}$ and such that $U \cap C_{0} \cap Sing (X_{0})
= \emptyset$, let $B_{1}(U)$ be the subvariety of $t \in B$
parametrizing projectively Cohen-Macaulay surfaces such that the
hypersurface $X_{t}$ contains $U \cap C_{0}$. Now let $B(U)$ be an
irreducible component of $B_{1}(U)$ containing 
$0$ and the points in $B_{1}(U)$ parametrizing projectively 
Cohen-Macaulay surfaces whose matrix is of type 
$$\left[\matrix{A_{11} & \ldots & A_{1m} & UG_{1} \cr \vdots & 
\null & \null & \vdots \cr A_{m1} & \ldots & A_{mm} & UG_{m}
\cr}\right].$$   
\noindent Let $T(U)$ be the tangent space to $B(U)$ at $0$ and 
$\widetilde {\cal C}(U), {\cal X}(U)$ the corresponding families of
curves and threefolds.  We will prove that the characteristic map of
$\widetilde{\cal C}(U)$ has rank  3 for $U$ general.
To prove this using Lemma (3.1), we will show first that if $U, U'$ 
are monomials which differ only in degree one, i.e. $V = $ 
l.c.m.$(U,U')$ has degree $u + 1$, then the characteristic map on
$B(U) \cap B(U')$ has rank at least 2. In fact by Proposition (2.1)
it is  enough to show that the corresponding curves $C_{t}, t \in
B(U) \cap B(U')$, fill up a variety 
$\overline{\Sigma} \subset \I^4$ of dimension at least 3. Suppose 
to the contrary $\hbox{dim} \overline{\Sigma} \leq 2$. As we can
assume that $V$ is transversal to $C_{0}$, it is necessarily
$\hbox{dim} \overline{\Sigma} 
\cap V \leq 1$, since $C_{0} \subset \overline{\Sigma}$. What we 
need follows then by the

\noindent {\it Claim (3.6).} {\sl For all monomials $U, U', V$ as 
above and for every variety $\Sigma \subset \I^4$ such that
$\hbox{dim} \Sigma \leq 2$  and $\hbox{dim} \Sigma \cap V \leq 1$,
there exists $t \in B(U) \cap B(U')$  such that $\hbox{dim} S_{t}
\cap \Sigma = 0$.}

\noindent {\it Proof of Claim (3.6)}: We choose $t \in B(U) \cap 
B(U')$ so that the matrix of $S_{t}$ is a general one of type 
$$\left[\matrix{A_{11} & \ldots & A_{1m} & VH_{1} \cr \vdots &
\null & \null & \vdots \cr A_{m1} & \ldots & A_{mm} & VH_{m} 
\cr}\right]$$
\noindent that is, the polynomials $A_{ij}$ and $H_{i}$ (of degree 
$u_{i,m+1} - u -1$) are general. Let $T$ be the projectively
Cohen-Macaulay surface  defined by the vanishing of the $(m-1)
\times (m-1)$ minors of the matrix
$[A_{ij}, 1 \leq i \leq m, 2 \leq j \leq m]$ and $F_{1}, F_{m+1}$ 
the two distinguished generators  of the ideal of $S_{t}$. Clearly
$\hbox{dim} T \cap \Sigma = 0, \hbox{dim}  F_{m+1} \cap \Sigma = 1$
and $F_{m+1}$ does not  contain any component of $V \cap \Sigma$, as
the entries $A_{ij}$ are general. Now suppose that $S_{t} \cap
\Sigma$ has a component $D$ of dimension at least one and let $x 
\in D$ be a general point. Then $x \not\in V$. As 
$D \not\subset T$, there is an $(m-1) \times (m-1)$ minor of the 
matrix of $T$ not vanishing on $x$; hence $x \not\in F_{1}$, by a
general choice of the $H_{i}$.  
\noindent This contradiction proves Claim (3.6). 

\noindent To finish the proof of the theorem just proceed in analogy
 with the previous proofs. Since $T_{B}$ is generated by the tangent
vectors to the varieties $B(U)$, with $U$ monomial and the
characteristic map of $\widetilde{\cal C}$ has rank 3 at a general
point, we get by Lemma (3.1) that the characteristic map of
$\widetilde{\cal C}(U)$ has also rank 3. Now look at  the focal
locus of $\widetilde{\cal C}(U)$: it contains the inverse image of
$C_{0} \cap U$, as they are fixed points of the family of 
hypersurfaces ${\cal X}(U)$; hence  by Propositions (2.3) and
(2.4), 
$$2g-2\ge u \hbox{deg} C_{0} - 5\hbox{deg} C_{0}. \ \ \ \qed$$ 
\noindent {\bf (3.7)}{\it Remark.} It is clear that if we drop the 
hypothesis $u_{m,m+1} \geq 8$ in Theorem (1.4) there can be rational
or elliptic curves on a general projectively Cohen-Macaulay surface
$S$ in $\I^4$. For example  the Castelnuovo surface contains
elliptic curves.

\noindent {\bf (3.8)}{\it Remark.} The Picard group of the surfaces 
considered in our theorems is, in many cases, particularly simple. 
In the case of Theorem (1.2) it is just generated by the hyperplane
bundle, by the Noether-Lefschetz theorem. In the hypothesis of
Theorem (1.3) it  follows by [Lo, Corollary II.3.8] that $S$ has
Picard group generated by the hyperplane bundle $H$ and by the line
bundle associated to $D$. For the general projectively
Cohen-Macaulay surface $S \subset \I^4$ the situation is quite
different. If $u_{ij} > 0$ for every $i,j$, the Picard group is
generated by the hyperplane bundle $H$ and by the canonical bundle
of $S$ [Lo, Theorem III.4.2], while in general it can have large
rank. Despite of the simplicity  of the Picard groups in many cases,
there does not seem to be a way to make use of it, as the theorems
are concerned with the possible singularities of the curves. As a
matter of fact our proofs are independent of the knowledge on the
Picard group.

\vskip1cm 

\centerline{\bf REFERENCES}

\baselineskip 12pt

\vskip.3cm
\item{[B]} Brody,R.:\ Compact manifolds in hyperbolicity.\ {\it 
Trans.\ Amer.\ Math.\ Soc.\ \bf 235}, (1978) 213-219.  
\vskip.16cm
\item{[BG]} Brody,R., Green,M.:\ A family of smooth hyperbolic 
hypersurfaces in $\I^3$.\ {\it Duke Math.\ J.\ \bf 44}, (1977)
873-874.   
\vskip.16cm
\item{[CC]} Chiantini,L., Ciliberto,C.:\ A few remarks on the 
lifting problem.\ {\it Asterisque\ \bf 218}, (1993) 95-109.
\vskip.16cm 
\item{[Ch]} Chang,M.C.:\ A filtered Bertini type theorem.\ 
{\it J.\ Reine Angew.\ Math.\ \bf 397}, (1989) 214-219.   
\vskip.16cm
\item{[Cl]} Clemens,H.:\ Curves on generic hypersurfaces.\ {\it 
Ann.\ Sci.\ \'Ecole Norm.\ Sup.\ \bf 19}, (1986) 629-636.   
\vskip.16cm 
\item{[CLR]} Chiantini,L., Lopez,A.F., Ran, Z.:\ Subvarieties of generic
hypersurfaces in any variety.\ {\it Preprint alg-geom AG/9901083.}
\vskip.16cm
\item{[CS]} Ciliberto,C., Sernesi,E.:\ Singularities of the theta 
divisor and congruences of planes.\ {\it J.\ Alg.\ Geom.\ \bf 1},
(1992) 231-250. 
\vskip.16cm 
\item{[D]} Demailly,J.P.:\ Algebraic criteria for Kobayashi 
hyperbolic projective varieties and jet differentials.\ {\it Notes 
of the AMS Summer Institute, Santa Cruz 1995.} 
\vskip.16cm 
\item{[DEG]} Demailly,J.P., El Goul,J.:\ Hyperbolicity of generic
surfaces of high degree in projective 3-space.\ {\it Preprint 
alg-geom AG/9804129.} 
\vskip.16cm 
\item{[E1]} Ein,L.:\ Subvarieties of generic complete 
intersections.\ {\it Invent.\ Math.\ \bf 94}, (1988) 163-169.
\vskip.16cm
\item{[E2]} Ein,L.:\ Subvarieties of generic complete intersections 
II.\ {\it  Math.\ Ann.\ \bf 289}, (1991) 465-471.
\vskip.16cm
\item{[EG]} El Goul,J.:\ Algebraic families of smooth hyperbolic 
surfaces of low degree in $\I^3_{\hbox{\msq C}}$.\ {\it Manuscripta
Math.\ \bf 90}, (1996) 521-532.
\vskip.16cm 
\item{[K]} Kobayashi,S.:\ Hyperbolic manifolds and holomorphic 
mappings.\ {\it Pure and Applied Mathematics \bf 2}.\ Marcel Dekker,
New York\ 1970.    
\vskip.16cm
\item{[La]} Lang,S.:\ Hyperbolic and Diophantine analysis.\ {\it 
Bull.\ Amer.\ Math.\ Soc.\ \bf 14}, (1986) 159-205.
\vskip.16cm
\item{[Lo]} Lopez,A.F.:\ Noether-Lefschetz theory and the Picard 
group of projective surfaces.\ {\it Mem.\ Amer.\ Math.\ Soc.\ \bf
89}, vol. 438 (1991).
\vskip.16cm  
\item{[N]} Nadel,A.:\ Hyperbolic surfaces in $\I^3$.\ {\it Duke 
Math.\ J.\ \bf 58}, (1989) 749-771.  
\vskip.16cm
\item{[S]} Segre,C.:\ Sui fuochi di $2^o$ ordine dei sistemi 
infiniti di piani, e sulle curve iperspaziali con una doppia
infinit\`a di piani plurisecanti.\ {\it Atti R.\ Accad.\ Lincei \bf
30}, vol. 5, (1921) 67-71.
\vskip.16cm
\item{[V]} Voisin,C.:\ On a conjecture of Clemens on rational 
curves on hypersurfaces.\ {\it J.\ Differential Geom.\ \bf 44},
(1996) 200-213. \ {\it J.\ Differential Geom.\ \bf 49}, 1998, 601-611.
\vskip.16cm
\item{[X1]} Xu,G.:\ Subvarieties of general hypersurfaces in 
projective space.\ {\it J.\ Differential Geom.\ \bf 39}, (1994)
139-172.    
\vskip.16cm
\item{[X2]} Xu,G.:\ Divisors on hypersurfaces.\ {\it Math.\ Z.\ 
\bf 219}, (1995) 581-589.
\vskip.16cm
\item{[X3]} Xu,G.:\ Divisors on general complete intersections in 
projective space.\ {\it Trans.\ Amer.\ Math.\ Soc.\ \bf 348}, (1996)
2725-2736.

\vskip.5cm

\baselineskip=12pt
\ \ \ \ \ \ \ LUCA CHIANTINI \ \ \ \ \ \ \ \ \ \ \ \ \ \ \ \ \ \ \ \ \ 
\ \ \ \ \ \ ANGELO FELICE LOPEZ

\ \ Dipartimento di Matematica \ \ \ \ \ \ \ \ \ \ \ \ \ \ \ \ \ \ 
\ \ 
Dipartimento di Matematica

\ \ \ \ \ \ \ Universit\`a di Siena \ \ \ \ \ \ \ \ \ \ \ \ \ \ \ 
\ \ \ \ \ \ \ \ \ 
\ \ \ \ \ \ \ Universit\`a di Roma Tre

\ \ \ \ \ \ Via del Capitano 15 \ \ \ \ \ \ \ \ \ \ \ \ \ \ \ \ \ 
\ \ \ \ \ \ \ \ \ Largo San Leonardo Murialdo 1

\ \ \ \ \ \ \ \ 53100 Siena  Italy \ \ \ \ \ \ \ \ \ \ \ \ \ \ \ \ 
\ \ \ \ \
\ \ \ \ \ \ \ \ \ \ \ \ \ \ \ \ 00146 Roma Italy 

\ \ e-mail {\tt{chiantini@unisi.it}}  \ \ \ \ \ \ \ \ \ \ \ \ 
\ \ e-mail {\tt{lopez@matrm3.mat.uniroma3.it}}  

\end